\documentclass[11pt]{article}%
\usepackage{amsmath}
\usepackage{amsfonts}
\usepackage{amssymb}
\usepackage{graphicx}
\usepackage{theorem}
\usepackage{color}%
\providecommand{\U}[1]{\protect\rule{.1in}{.1in}}
\newtheorem{thm}{Theorem}[section]
\newtheorem{lm}[thm]{Lemma}
\newtheorem{pr}[thm]{Proposition}
\newtheorem{df}[thm]{Definition}
\newtheorem{rmk}[thm]{Remark}
\newtheorem{cor}[thm]{Corollary}
{\theorembodyfont{\upshape}
\newtheorem{examp}[thm]{Example}
}
\numberwithin{equation}{section} 
\begin{document}

\title{Title: }

\begin{center}
\bigskip

\vspace*{1.3cm}

\textbf{SUBDIFFERENTIAL CALCULUS AND IDEAL\ SOLUTIONS\ FOR SET
OPTIMIZATION\ PROBLEMS}

\bigskip

by

\bigskip

Marius\ DUREA\footnote{{\small Faculty of Mathematics, \textquotedblleft
Alexandru Ioan Cuza\textquotedblright\ University, 700506--Ia\c{s}i, Romania
and \textquotedblleft Octav Mayer\textquotedblright\ Institute of Mathematics,
Ia\c{s}i Branch of Romanian Academy, 700505--Ia\c{s}i, Romania{; e-mail:
\texttt{durea@uaic.ro}}}} and Elena-Andreea FLOREA\footnote{{\small Faculty of
Mathematics, \textquotedblleft Alexandru Ioan Cuza\textquotedblright%
\ University, 700506--Ia\c{s}i, Romania and \textquotedblleft Octav
Mayer\textquotedblright\ Institute of Mathematics, Ia\c{s}i Branch of Romanian
Academy, 700505--Ia\c{s}i, Romania{; e-mail: {\small \texttt{a}}%
\texttt{ndreea.florea@uaic.ro}}}}
\end{center}

\bigskip

\noindent{\small {\textbf{Abstract:}} We explore the possibility to derive
basic calculus rules for some subdifferential constructions associated to
set-valued maps between normed vector spaces. Then, we use these results in
order to write optimality conditions for a special kind of solutions for set
optimization problems. }

\bigskip


\bigskip


\begin{center}

\end{center}

\section{Introduction}

The primary aim of this work is to develop some ideas concerning the
subdifferential calculus rules for generalized differentiation objects
associated to set-valued maps (or multifunctions) which were introduced and
only briefly pointed out in the paper \cite{DS23-JOGO}. The generalized
subgradients under study in this paper were initially designed to deal with
set optimization problems and therefore are constructed on the basis of
epigraphical set-valued maps, where the epigraphs are defined by means of an
ordering cone in the output space. Actually, one of the main tools we
systematically employ for deriving our results is the scalarization of the
underlying set-valued maps by elements in the dual cone associated to the
ordering cone and in doing so, and further relying on some topological
conditions, we reduce, in some relevant cases, the calculus for the
subgradients of the set-valued maps with vectorial values to the subgradients
of real-valued functions, a case for which, of course, many calculus rules are
available (see \cite{Morduk2006}). Besides the Fr\'{e}chet subdifferential of
set-valued maps introduced in \cite{DS23-JOGO}, we define here a corresponding
limiting (or Mordukhovich) subdifferential and for both types of subgradients
we study the (generalized) convex case and sum rules. Needless to say, the
calculus we derive here covers known results from the case of nonsmooth
real-valued functions. Moreover, a difference rule for Fr\'{e}chet
subdifferential is derived in a less general situation but this case is,
however, general enough to ensure meaningful necessary optimality conditions
for a newly introduced type of solution for set optimization problems. We
notice that the investigation of the main results is based on several
statements that could be of some importance for their own: we mention here a
"conic" variant of R\aa dstr\"{o}m cancellation law and a penalization
procedure. Furthermore, the concept of (sequential) cone compactness of a set
seems to be an instrument in various situations that naturally appear in our study.

The paper is organized as follows. The second section collects some
preliminaries concerning scalarization by positive linear continuous
functionals and useful (semi)continuity properties of set-valued maps.

The main section is the third one. Firstly, we present several calculus rules
for already introduced Fr\'{e}chet subdifferential of set-valued maps between
normed vector spaces. As we already mentioned, on the output space we consider
a closed convex pointed cone and we start our development of calculus rules by
the case of convex multifunctions (where the convexity is defined by means of
the ordering cone). In this basic situation, it becomes apparent that the
study requires an adapted version of R\aa dstr\"{o}m cancellation law
involving the presence of the ordering cone. Moreover, several concepts
related to sets that are used thoroughly afterwards are recalled, the most
important one being the sequential compactness with respect to a cone
introduced and studied in \cite{DF2022}. Then we present a sum rule in the
convex case. Secondly, on the basis of the Fr\'{e}chet subdifferential, we
introduce in the same manner (i.e., by passing to the limit procedure) as in
the case of real-valued functions the limiting (or Mordukhovich)
subdifferential for which we get formulas in the convex case and for the sum.
Thirdly, we are interested in some calculus rules concerning the
subdifferentials of some special (multi)functions and, in particular, we
derive an estimation formula for the Fr\'{e}chet subdifferential of a
difference between a general set-valued map and a particular function.

The last section applies some of the calculus rules developed in the previous
section in order to get necessary optimality conditions for a special type of
solutions for set optimization problems, which we call ideal solutions. In
order to arrive at the necessary optimality conditions we want to display,
several technical results, among which we mention a penalization procedure,
are derived.

\bigskip

We briefly present the notation we use in this work. Let $X,Y$ be normed
spaces over the real field $\mathbb{R}$. The topological dual of $X$ is
$X^{\ast},$ the norm will be denoted $\left\Vert \cdot\right\Vert ,$ and
$B\left(  X,Y\right)  $ is the normed vector space of linear bounded operators
from $X$ to $Y.$ For $x\in X$ and $\varepsilon>0,$ we put $B\left(
x,\varepsilon\right)  $ for the open ball centered at $x$ with the radius
$\varepsilon.$ The symbols $S_{X}$ and $D_{X}$ stand for the unit sphere and
for the closed unit ball, respectively. If $A\subset X$ is a nonempty set,
then we denote by $\operatorname{cl}A,$ $\operatorname*{int}A,$ $A^{\prime},$
$\operatorname*{conv}A,$ $d\left(  \cdot,A\right)  $ the topological closure,
the topological interior, the set of accumulation points, the convex hull, and
the associated distance function, respectively. If $A,B$ are nonempty subsets
of $X$, the excess from $A$ to $B$ is $e(A,B)=\sup_{x\in A}d(x,B).$

\section{Preliminaries}

In general, in the following, $K$ is supposed to be a pointed closed convex
cone in $Y$ and its positive dual cone is denoted by $K^{+}.$ We consider a
set-valued map $F:X\rightrightarrows Y$ and we associate with it the
epigraphical set-valued map $\operatorname*{Epi}F:X\rightrightarrows Y$ given
by $\operatorname*{Epi}F\left(  x\right)  =F\left(  x\right)  +K.$ As a
standing assumption, we assume that all the values of $F$ that are involved in
the discussion are nonempty sets. The following concepts were introduced in
\cite{DS23-JOGO}.

\begin{df}
\label{subdiff}Let $F:X\rightrightarrows Y,$ and $\overline{x}\in X.$ The
Fr\'{e}chet subdifferential of $F$ at $\overline{x}$ is%
\begin{equation}
\widehat{\partial}F\left(  \overline{x}\right)  =\left\{  T\in B\left(
X,Y\right)  \mid\lim_{x\rightarrow\overline{x}}\frac{e\left(
\operatorname*{Epi}F\left(  x\right)  ,\operatorname*{Epi}F\left(
\overline{x}\right)  +T\left(  x-\overline{x}\right)  \right)  }{\left\Vert
x-\overline{x}\right\Vert }=0\right\}  . \label{Fr}%
\end{equation}
Equivalently, $T\in\widehat{\partial}F\left(  \overline{x}\right)  $ iff $T\in
B\left(  X,Y\right)  $ and%
\begin{equation}
\forall\varepsilon>0,\exists\delta>0,\forall x\in B\left(  \overline{x}%
,\delta\right)  :\operatorname*{Epi}F\left(  x\right)  \subset
\operatorname*{Epi}F\left(  \overline{x}\right)  +T\left(  x-\overline
{x}\right)  +\varepsilon\left\Vert x-\overline{x}\right\Vert D_{Y}.
\label{Fr2}%
\end{equation}

Similarly, we can define the upper subdifferential of $F$ at $\overline{x}$ as
follows%
\begin{equation}
\widehat{\partial}^{+}F\left(  \overline{x}\right)  =\left\{  T\in B\left(
X,Y\right)  \mid\lim_{x\rightarrow\overline{x}}\frac{e\left(
\operatorname*{Epi}F\left(  \overline{x}\right)  +T\left(  x-\overline
{x}\right)  ,\operatorname*{Epi}F\left(  x\right)  \right)  }{\left\Vert
x-\overline{x}\right\Vert }=0\right\}  . \label{uFr}%
\end{equation}

\end{df}

As was remarked in the same work, these notions generalize the respective
subdifferentials from the case of scalar functions. Indeed, if $f:X\rightarrow
\mathbb{R}$ is a function, then relation (\ref{Fr2}) can be equivalently
written as%
\[
\forall\varepsilon>0,\exists\delta>0,\forall x\in B\left(  \overline{x}%
,\delta\right)  :f\left(  x\right)  -f\left(  \overline{x}\right)  -x^{\ast
}\left(  x-\overline{x}\right)  \geq-\varepsilon\left\Vert x-\overline
{x}\right\Vert ,
\]
i.e., $x^{\ast}\in\widehat{\partial}f\left(  \overline{x}\right)  ,$ where
$\widehat{\partial}f\left(  \overline{x}\right)  $ denotes the usual
Fr\'{e}chet subdifferential of $f$ at $\overline{x}$ (see \cite{Morduk2006}).
A similar comment holds for $\widehat{\partial}^{+}f\left(  \overline
{x}\right)  .$

\bigskip

The next result is based on \cite[Lemma 2.5]{DS23-JOGO}. For the completeness,
we give a direct proof and we point out as well a missing assumption in the
mentioned result.

\begin{pr}
\label{rez_incl}Let $K\subset Y$ be a pointed closed convex cone, and $A,B$
nonempty subsets of $Y.$

(i) If $A\subset B+K,$ then
\begin{equation}
y^{\ast}\left(  A\right)  +[0,+\infty)\subset y^{\ast}\left(  B\right)
+[0,+\infty),\text{ }\forall y^{\ast}\in K^{+}\setminus\left\{  0\right\}  .
\label{rez_incl1}%
\end{equation}
The converse implication holds provided $B+K$ is convex and closed.

(ii) If $\operatorname*{int}K\neq\emptyset$ and $A\subset
B+\operatorname*{int}K,$ then
\[
y^{\ast}\left(  A\right)  \subset y^{\ast}\left(  B\right)  +(0,+\infty
),\text{ }\forall y^{\ast}\in K^{+}\setminus\left\{  0\right\}  .
\]
The converse implication holds provided $B+K$ is convex.
\end{pr}

\noindent\textbf{Proof. }If $A\subset B+K,$ then for all $a\in A,\,$there is
$b\in B\ $such that $a\in b+K.$ Therefore, for all $y^{\ast}\in K^{+}%
\setminus\left\{  0\right\}  $ one has $y^{\ast}\left(  b\right)  \leq
y^{\ast}\left(  a\right)  ,$ meaning that
\[
y^{\ast}\left(  a\right)  \in y^{\ast}\left(  b\right)  +[0,\infty)\subset
y^{\ast}\left(  B\right)  +[0,\infty).
\]
Since $a$ was arbitrarily chosen in $A,$ we get $y^{\ast}\left(  A\right)
+[0,+\infty)\subset y^{\ast}\left(  B\right)  +[0,+\infty)$ for all $y^{\ast
}\in K^{+}\setminus\left\{  0\right\}  .$

Suppose that $B+K$ is convex and closed. Suppose, by way of contradiction,
that the converse does not hold. Then there is $a\in A$ such that $a\notin
B+K.$ The assumptions we made allow us to strongly separate the point $a$ from
the set $B+K,$ so one gets $y^{\ast}\in Y^{\ast}\setminus\left\{  0\right\}  $
and $\alpha\in\mathbb{R}$ such that%
\[
y^{\ast}\left(  a\right)  <\alpha<y^{\ast}\left(  b+k\right)  ,\text{ }\forall
b\in B,\text{ }\forall k\in K.
\]
Clearly, this forces $y^{\ast}\in K^{+}\setminus\left\{  0\right\}  $ and
\[
y^{\ast}\left(  a\right)  <\alpha<y^{\ast}\left(  b\right)  ,\text{ }\forall
b\in B.
\]
But, by assumption, there is $b_{a}\in B$ such that $y^{\ast}\left(
b_{a}\right)  \leq y^{\ast}\left(  a\right)  $ and we reach a contradiction.

(ii) Since $y^{\ast}\left(  k\right)  >0$ for all $y^{\ast}\in K^{+}%
\setminus\left\{  0\right\}  $ and $k\in\operatorname*{int}K,$ we easily get
the direct implication. For the converse, suppose that $B+K$ is convex and
there is $a\in A$ such that $a\notin B+\operatorname*{int}K.$ Therefore, one
gets $y^{\ast}\in Y^{\ast}\setminus\left\{  0\right\}  $ such that%
\[
y^{\ast}\left(  a\right)  <y^{\ast}\left(  b+k\right)  ,\text{ }\forall b\in
B,\text{ }\forall k\in\operatorname*{int}K.
\]
Again, this forces $y^{\ast}\in K^{+}$ and
\[
y^{\ast}\left(  a\right)  \leq y^{\ast}\left(  b\right)  ,\text{ }\forall b\in
B.
\]
But, by assumption, there is $b_{a}\in B$ such that $y^{\ast}\left(
b_{a}\right)  <y^{\ast}\left(  a\right)  $ and this contradicts the above
inequality.\hfill$\square$

\begin{examp}
Let $Y=\mathbb{R}^{2},$ $A=\left\{  \left(  1,1\right)  \right\}  ,$
$B=\left\{  \left(  2,0\right)  ,\left(  0,2\right)  \right\}  $ and
$K=\mathbb{R}_{+}^{2}.$ Then for $y^{\ast}=\left(  y_{1},y_{2}\right)  \in
K^{+}\setminus\left\{  0\right\}  =\mathbb{R}_{+}^{2}\setminus\left\{
0\right\}  ,$
\[
y^{\ast}\left(  1,1\right)  =y_{1}+y_{2}\geq\left\{
\begin{array}
[c]{c}%
2y_{1},\text{ if }y_{2}\geq y_{1}\\
2y_{2},\text{ if }y_{1}>y_{2}.
\end{array}
\right.
\]
Since $2y_{1},$ $2y_{2}\in y^{\ast}\left(  B\right)  ,$ we see that the
converse in Proposition \ref{rez_incl} (i) does not hold without additional assumptions.
\end{examp}

\begin{rmk}
Of course, for the converses in $(i)$ and $(ii)$ of the above result, one can
take $y^{\ast}\in S_{Y^{\ast}}\cap K^{+}$ or only $y^{\ast}$ in a set of
generators of $K^{+}$ without $0.$
\end{rmk}

Then, we recall the notions of lower and upper continuity of a set-valued map
(see \cite[Definition 2.5.1]{GoeRiaTamZal2003}) and a generalized Lipschitz
condition with respect to $K$ of a set-valued map (see \cite{Ye}, \cite{BDS}).

\begin{df}
Let $F:X\rightrightarrows Y$ be a set-valued map with nonempty values and
$x\in X$. Then $F$ is said to be:

(i) lower continuous at $x$ (in short, l.c.) if for each open set $V\subset Y$
with $F\left(  x\right)  \cap V\neq\emptyset$, there is a neighborhood
$U\subset X$ of $x$ such that for each $u\in U\ $one has $F\left(  u\right)
\cap V\neq\emptyset$;

(ii) upper continuous at $x$ (in short, u.c.) if for each open set $V\subset
Y$ with $F\left(  x\right)  \subset V$, there is a neighborhood $U\subset X$
of $x$ such that for each $u\in U\ $one has $F\left(  u\right)  \subset V$;

(iii) continuous at $x$ if it is both l.c. and u.c. at $x;$

(iv) l.c. (u.c., continuous) on a nonempty set $A\subset X$ if it is l.c.
(u.c., continuous) at all $x\in A$;

(v) l.c. (u.c., continuous) around $x$ if there exists $\varepsilon>0$ such
that $F$ is l.c. (u.c., continuous) on $B\left(  x,\varepsilon\right)  .$
\end{df}

\begin{df}
\label{def_K_lip}Let $F:X\rightrightarrows Y$ be\ a set-valued map and
$\overline{x}\in X$. One says that $F$ is $K-$Lipschitz around $\overline{x}$
if there exist a neighborhood $U$ of $\overline{x}$, a constant $L>0$ and an
element $e\in K\setminus\{0\}$ such that for every $x,u\in U,$%
\[
F(x)+L\left\Vert x-u\right\Vert e\subset F(u)+K.
\]

\end{df}

\section{Subdifferential calculus}

We intend to derive some subdifferential calculus rules for the Fr\'{e}chet
subdifferential and then to introduce and study an associated limiting
(Mordukhovich) subdifferential.

Let us first remark that if $F$ is a set-valued map from $X$ to $\mathbb{R}$
and has closed, bounded from below values, then taking $f:X\rightarrow
\mathbb{R}$ given by $f\left(  x\right)  =\min F\left(  x\right)  ,$ one has,
with the natural choice $K=[0,\infty)$, for every $\overline{x}\in X,$ that
$t\in\widehat{\partial}F\left(  \overline{x}\right)  $ if and only if
$t\in\widehat{\partial}f\left(  \overline{x}\right)  ,$ and similarly for
$\widehat{\partial}^{+}.$ Notice as well that for this kind of set-valued map,
the boundedness from below is required in order to avoid trivial situations in
terms of $\widehat{\partial}$ and $\widehat{\partial}^{+}.$

From this point on we concentrate on $\widehat{\partial}.$ The first result we
derive considers a scalarization of a general set-valued map by elements in
the positive dual cone in order to get the situation described above (that is,
the case of set-valued maps taking as values nonempty subsets of $\mathbb{R}$).

\begin{pr}
\label{pr_scal_subd}Let $F:X\rightrightarrows Y$ be a set-valued map,
$\overline{x}\in X,$ and $T\in\widehat{\partial}F\left(  \overline{x}\right)
$. Then for all $y^{\ast}\in K^{+}\setminus\left\{  0\right\}  $ one has
$y^{\ast}\circ T\in\widehat{\partial}\left(  y^{\ast}\circ F\right)  \left(
\overline{x}\right)  .$ If $F\left(  \overline{x}\right)  +K$ is convex and
closed, and the cone $K$ is finitely generated, then the converse holds.
\end{pr}

\noindent\textbf{Proof. }If $T\in\widehat{\partial}F\left(  \overline
{x}\right)  ,$ then for all $\varepsilon>0,$ there is $\delta>0$ such that for
all$\ x\in B\left(  \overline{x},\delta\right)  $ one has%
\[
F\left(  x\right)  \subset F\left(  \overline{x}\right)  +T\left(
x-\overline{x}\right)  +\varepsilon\left\Vert x-\overline{x}\right\Vert
D_{Y}+K.
\]
By Proposition \ref{rez_incl}, for all $y^{\ast}\in K^{+}\setminus\left\{
0\right\}  ,$%
\begin{equation}
\left(  y^{\ast}\circ F\right)  \left(  x\right)  \in\left(  y^{\ast}\circ
F\right)  \left(  \overline{x}\right)  +\left(  y^{\ast}\circ T\right)
\left(  x-\overline{x}\right)  +\varepsilon\left\Vert x-\overline
{x}\right\Vert y^{\ast}\left(  D_{Y}\right)  +[0,\infty),\text{ }\forall x\in
B\left(  \overline{x},\delta\right)  , \label{rel_scal}%
\end{equation}
which implies%
\[
\left(  y^{\ast}\circ F\right)  \left(  x\right)  \in\left(  y^{\ast}\circ
F\right)  \left(  \overline{x}\right)  +\left(  y^{\ast}\circ T\right)
\left(  x-\overline{x}\right)  -\varepsilon\left\Vert x-\overline
{x}\right\Vert \left\Vert y^{\ast}\right\Vert +[0,\infty),\text{ }\forall x\in
B\left(  \overline{x},\delta\right)  .
\]
Of course, this is enough to conclude that $y^{\ast}\circ T\in\widehat
{\partial}\left(  y^{\ast}\circ F\right)  \left(  \overline{x}\right)  .$

For the converse, if $K$ is finitely generated, then $K^{+}$ is finitely
generated (see, for instance, \cite[Section 2.3]{Ber}) and since in inclusion
(\ref{rez_incl1}) in Proposition \ref{rez_incl} we can take only elements from
a set of generators of $K^{+}$ and this is finite, when we write $y^{\ast
}\circ T\in\widehat{\partial}\left(  y^{\ast}\circ F\right)  \left(
\overline{x}\right)  $, for all $\varepsilon>0$ we get the same $\delta>0$
working for all $y^{\ast}$ in that set of generators for which (\ref{rel_scal}%
) holds and we can apply Proposition \ref{rez_incl}.\hfill$\square$

\bigskip

Our primary concern is to see that, as in the classical case, the
subdifferentials take a special form under convexity assumptions, so we
consider now the convex case. Here, the convexity is understood in the
generalized sense described in the definition below.

\begin{df}
Let $A\subset X$ be a nonempty convex set and $G:A\rightrightarrows Y$ be\ a
set-valued map.

(i) One says that $G\ $is upper $K-$convex if for every $x,y\in A$ and every
$\lambda\in\left(  0,1\right)  $ one has%
\[
\lambda G\left(  x\right)  +\left(  1-\lambda\right)  G\left(  y\right)
\subset G\left(  \lambda x+\left(  1-\lambda\right)  y\right)  +K.
\]

(ii) Suppose that $K$ is solid (that is, $\operatorname*{int}K\neq\emptyset$).
One says that $G$ is strictly upper $K-$convex if for every $x,y\in A$ and
every $\lambda\in\left(  0,1\right)  $ one has%
\[
\lambda G\left(  x\right)  +\left(  1-\lambda\right)  G\left(  y\right)
+\operatorname*{int}K\subset G\left(  \lambda x+\left(  1-\lambda\right)
y\right)  +\operatorname*{int}K.
\]

\end{df}

Inspired by \cite[Theorem 2.2]{Sch} we obtain the next result which is a
variant of R\aa dstr\"{o}m cancellation lemma adapted to our setting. We use
the following notion: a nonempty set $A\subset Y$ is called $K-$bounded if
there is a bounded set $M\subset Y$ such that $A\subset M+K.$

\begin{lm}
\label{lm_Rad}Suppose that $A,B,C\subset Y\ $are nonempty sets such that $C$
is $K-$bounded and
\[
A+C\subset\operatorname{cl}\left(  C+B+K\right)  .
\]
Then
\[
A\subset\operatorname{cl}\operatorname*{conv}\left(  B+K\right)  .
\]

\end{lm}

\noindent\textbf{Proof. }The proof is similar to one of the proofs of
R\aa dstr\"{o}m cancellation lemma (see, for instance, \cite[Theorem 2.2]%
{Sch}). Firstly, observe that it is enough to take $A$ as a singleton,
$A=\left\{  a\right\}  $ with $a\in Y.$ Then, since $\left\{  a\right\}
+C\subset\operatorname{cl}\left(  C+B+K\right)  ,$ one has $C\subset
\operatorname{cl}\left(  C+B+K\right)  -a=\operatorname{cl}\left(  C+\left(
B-a\right)  +K\right)  $ and as the following equalities hold
\[
\operatorname*{cl}\operatorname*{conv}\left(  \left(  B-a\right)  +K\right)
=\operatorname*{cl}\left(  \operatorname*{conv}\left(  B+K\right)  -a\right)
=\operatorname{cl}\operatorname*{conv}\left(  B+K\right)  -a,
\]
we can even consider $A=\left\{  0\right\}  .$

Next, we know that%
\[
C\subset\operatorname{cl}\left(  C+B+K\right)
\]
and we have to show that $0\in\operatorname{cl}\operatorname*{conv}\left(
B+K\right)  .$ Denote by $M$ a bounded set that satisfies $C\subset M+K.$ Then
for all natural $n\geq1,$%
\begin{align*}
C &  \subset\operatorname{cl}\left(  C+B+K\right)  \subset\operatorname{cl}%
\left(  \operatorname{cl}\left(  C+B+K\right)  +B+K\right)  =\operatorname{cl}%
\left(  C+B+K+B+K\right)  \subset...\\
&  \subset\operatorname{cl}\left(  \underset{n\text{ times}}{\underbrace
{B+...+B}}+K+C\right)  \subset\operatorname{cl}\left(  \underset{n\text{
times}}{\underbrace{B+...+B}}+K+M\right)  .
\end{align*}
Take $c\in C.$ Then there is $c_{n}\in B\left(  c,n^{-1}\right)  $ and
$m_{n}\in M$ such that for all $n,$%
\[
\frac{1}{n}\left(  c_{n}-m_{n}\right)  \in\frac{1}{n}\left(  \underset{n\text{
times}}{\underbrace{\left(  B+K\right)  +...+\left(  B+K\right)  }}\right)
\subset\operatorname*{conv}\left(  B+K\right)  .
\]
Since $M$ is bounded, we get that $0\in\operatorname{cl}\operatorname*{conv}%
\left(  B+K\right)  $ and this is the conclusion.\hfill$\square$

\begin{rmk}
Of course, for $K=\left\{  0\right\}  $ the above results reduces to
\cite[Theorem 2.2]{Sch}.
\end{rmk}

Now we present the form of the Fr\'{e}chet subdifferential for a $K-$convex
set-valued map. We recall (see \cite{DF2022}) that a nonempty subset $A\subset
Y$ is called $K-$sequentially compact if for any sequence $\left(
a_{n}\right)  \subset A$ there is a sequence $\left(  c_{n}\right)  \subset K$
such that the sequence $\left(  a_{n}-c_{n}\right)  $ has a convergent
subsequence towards an element of $A$. It was proved in \cite{DF2022} that a
$K-$compact set (see \cite{Luc} for definition) is $K-$sequentially compact
and the converse holds provided $K\ $is separable. According to \cite{DF2022},
if $A$ is $K-$sequentially compact, then $A$ is also $K-$closed (that is,
$A+K$ is closed) and $K-$bounded. We also recall that $A$ is said to be
$K-$convex if $A+K$ is a convex set.

\begin{pr}
\label{pr_subd_F_convex}Suppose that $F$ is upper $K-$convex. Consider
$\overline{x}\in X$ and suppose that $F\left(  \overline{x}\right)  $ is
$K-$sequentially compact. Moreover, suppose that $D_{Y}$ is $K-$closed. Then
\[
\widehat{\partial}F\left(  \overline{x}\right)  =\left\{  T\in B\left(
X,Y\right)  \mid F\left(  x\right)  \subset F\left(  \overline{x}\right)
+T\left(  x-\overline{x}\right)  +K,\text{ }\forall x\in X\right\}  .
\]

\end{pr}

\noindent\textbf{Proof. }The fact that $\widehat{\partial}F\left(
\overline{x}\right)  $ includes the right-hand set is clear. Take
$T\in\widehat{\partial}F\left(  \overline{x}\right)  .$ Then, for all
$\varepsilon>0,$ there is $\delta_{\varepsilon}>0$ such that for all$\ x\in
B\left(  \overline{x},\delta_{\varepsilon}\right)  $ one has%
\[
F\left(  x\right)  \subset F\left(  \overline{x}\right)  +T\left(
x-\overline{x}\right)  +\varepsilon\left\Vert x-\overline{x}\right\Vert
D_{Y}+K.
\]
For all $x\in X,$ there is $\lambda_{\varepsilon,x}\in\left(  0,1\right)  $
such that $\left(  1-\lambda_{\varepsilon,x}\right)  \overline{x}%
+\lambda_{\varepsilon,x}x\in B\left(  \overline{x},\delta_{\varepsilon
}\right)  ,$ so we can write, taking into account the upper $K-$convexity of
$F$, that%
\begin{align*}
\left(  1-\lambda_{\varepsilon,x}\right)  F\left(  \overline{x}\right)
+\lambda_{\varepsilon,x}F\left(  x\right)   &  \subset F\left(  \left(
1-\lambda_{\varepsilon,x}\right)  \overline{x}+\lambda_{\varepsilon
,x}x\right)  +K\\
&  \subset F\left(  \overline{x}\right)  +\lambda_{\varepsilon,x}T\left(
x-\overline{x}\right)  +\varepsilon\lambda_{\varepsilon,x}\left\Vert
x-\overline{x}\right\Vert D_{Y}+K,
\end{align*}
whence%
\[
\left(  1-\lambda_{\varepsilon,x}\right)  F\left(  \overline{x}\right)
+\lambda_{\varepsilon,x}F\left(  \overline{x}\right)  +\lambda_{\varepsilon
,x}F\left(  x\right)  +K\subset F\left(  \overline{x}\right)  +\lambda
_{\varepsilon,x}F\left(  \overline{x}\right)  +\lambda_{\varepsilon,x}T\left(
x-\overline{x}\right)  +\varepsilon\lambda_{\varepsilon,x}\left\Vert
x-\overline{x}\right\Vert D_{Y}+K.
\]
But, the upper $K-$convexity of $F$ implies the convexity of the set $F\left(
\overline{x}\right)  +K,$ so%
\[
F\left(  \overline{x}\right)  +\lambda_{\varepsilon,x}F\left(  x\right)
+K\subset F\left(  \overline{x}\right)  +\lambda_{\varepsilon,x}F\left(
\overline{x}\right)  +\lambda_{\varepsilon,x}T\left(  x-\overline{x}\right)
+\varepsilon\lambda_{\varepsilon,x}\left\Vert x-\overline{x}\right\Vert
D_{Y}+K,
\]
whence%
\[
F\left(  \overline{x}\right)  +\lambda_{\varepsilon,x}F\left(  x\right)
\subset F\left(  \overline{x}\right)  +\lambda_{\varepsilon,x}F\left(
\overline{x}\right)  +\lambda_{\varepsilon,x}T\left(  x-\overline{x}\right)
+\varepsilon\lambda_{\varepsilon,x}\left\Vert x-\overline{x}\right\Vert
D_{Y}+K.
\]

Now, the set $F\left(  \overline{x}\right)  $ is $K-$bounded and the set
\[
\lambda_{\varepsilon,x}F\left(  \overline{x}\right)  +\lambda_{\varepsilon
,x}T\left(  x-\overline{x}\right)  +\varepsilon\lambda_{\varepsilon
,x}\left\Vert x-\overline{x}\right\Vert D_{Y}=\lambda_{\varepsilon,x}\left(
F\left(  \overline{x}\right)  +T\left(  x-\overline{x}\right)  +\varepsilon
\left\Vert x-\overline{x}\right\Vert D_{Y}\right)
\]
is $K-$closed (as a sum between a $K-$sequentially compact set and a
$K-$closed one: see \cite{DF2022}) and $K-$convex. Therefore, by Lemma
\ref{lm_Rad},%
\[
\lambda_{\varepsilon,x}F\left(  x\right)  \subset\lambda_{\varepsilon
,x}\left(  F\left(  \overline{x}\right)  +T\left(  x-\overline{x}\right)
+\varepsilon\left\Vert x-\overline{x}\right\Vert D_{Y}\right)  +K.
\]
We get that
\[
F\left(  x\right)  \subset F\left(  \overline{x}\right)  +T\left(
x-\overline{x}\right)  +\varepsilon\left\Vert x-\overline{x}\right\Vert
D_{Y}+K.
\]

For fixed $x$, this is true for all $\varepsilon$ and employing the
$K-$closedness of the set $F\left(  \overline{x}\right)  +T\left(
x-\overline{x}\right)  ,$ we finally deduce%
\[
F\left(  x\right)  \subset F\left(  \overline{x}\right)  +T\left(
x-\overline{x}\right)  +K,\text{ }\forall x\in X.
\]
Consequently, the equality is proven.\hfill$\square$

\begin{cor}
If $f:X\rightarrow Y$ is an upper $K-$convex function, $D_{Y}$ is $K-$closed,
and $\overline{x}\in X,$ then
\[
\widehat{\partial}f\left(  \overline{x}\right)  =\left\{  T\in B\left(
X,Y\right)  \mid f\left(  x\right)  \in f\left(  \overline{x}\right)
+T\left(  x-\overline{x}\right)  +K,\text{ }\forall x\in X\right\}  .
\]
In particular, if $Y=\mathbb{R}$ and $K=[0,\infty),$ this coincides with the
classical subdifferential of a convex function.
\end{cor}

\bigskip

In the previous results we use the assumption that $D_{Y}+K$ is a closed set.
Let us mention that this always happens when $Y$ is reflexive: indeed, in this
case $D_{Y}$ is weakly compact and $K$ is weakly closed, so the sum $D_{Y}+K$
is weakly closed, hence (norm-)closed. However, if $Y$ is not reflexive this
can fail to be true, as the next example shows.

\begin{examp}
Let $X=c_{0}$ with its usual (supremum) norm and define $f:c_{0}%
\rightarrow\mathbb{R},$
\[
f\left(  x\right)  =\sum_{n=1}^{\infty}\frac{1}{2^{n}}x_{n},\text{ }\forall
x=\left(  x_{n}\right)  _{n\geq1}\in c_{0}.
\]
It is easy to see that $f\in\left(  c_{0}\right)  ^{\ast},$ $\left\Vert
f\right\Vert =1$ and there is no $x\in D_{c_{0}}$ such that $f\left(
x\right)  =1,$ so $D_{c_{0}}\subset\left\{  x\in c_{0}\mid\left\vert f\left(
x\right)  \right\vert <1\right\}  .$

Consider the closed convex pointed cone
\[
A=\left\{  x=\left(  x_{n}\right)  _{n\geq1}\in c_{0}\mid x_{1}\leq0,\text{
}x_{n}\geq0,\text{ }\forall n\geq2\right\}  ,
\]
and define $K=A\cap\left\{  x\in c_{0}\mid f\left(  x\right)  \geq0\right\}
,$ which is also a closed convex pointed cone. One has that $D_{c_{0}%
}+K\subset\left\{  x\in c_{0}\mid f\left(  x\right)  >-1\right\}  .$ Consider
now the following sequences in $c_{0}:$ $\left(  x^{k}\right)  _{k\geq1},$
$\left(  y^{k}\right)  _{k\geq1}$ given by%
\[
x_{n}^{k}=\left\{
\begin{array}
[c]{l}%
-1,\text{ for }n\in\overline{1,k}\\
0,\text{ for }n>k
\end{array}
\right.
\]
and%
\[
y_{n}^{k}=\left\{
\begin{array}
[c]{l}%
-1,\text{ for }n=1\\
\left(  1+k^{-1}\right)  ,\text{ for }n\in\overline{2,k}\\
0,\text{ for }n>k.
\end{array}
\right.
\]
Then $\left(  x^{k}\right)  _{k}\subset D_{c_{0}}$ for all $k\geq1$ and
$\left(  y^{k}\right)  _{k}\subset K$ for all $k\geq3.$ Indeed, the former
inclusion is obvious, while for the latter we have the following computation:%
\begin{align*}
f\left(  y^{k}\right)   &  =-\frac{1}{2}+\sum_{n=2}^{k}\left(  1+\frac{1}%
{k}\right)  \frac{1}{2^{n}}=-\frac{1}{2}+\frac{1}{2}\left(  1+\frac{1}%
{k}\right)  \left(  1-\frac{1}{2^{k-1}}\right)  \\
&  =\frac{1}{2}\left(  \frac{1}{k}-\frac{1}{k}\frac{1}{2^{k-1}}-\frac
{1}{2^{k-1}}\right)  =\frac{1}{2}\frac{2^{k-1}-k-1}{k\cdot2^{k-1}}\geq0,\text{
}\forall k\geq3.
\end{align*}
So $\left(  x^{k}+y^{k}\right)  \subset D_{c_{0}}+K$ for $k\geq3.$ One has
that
\[
x_{n}^{k}+y_{n}^{k}=\left\{
\begin{array}
[c]{l}%
-2,\text{ for }n=1\\
k^{-1},\text{ for }n\in\overline{2,k}\\
0,\text{ for }n>k
\end{array}
\right.
\]
and, consequently, $x^{k}+y^{k}\rightarrow u=\left(  -2,0,0,...,0...\right)
.$ But, $f\left(  u\right)  =-1,$ whence $u\notin D_{c_{0}}+K,$ and this
proves that $D_{c_{0}}+K$ is not closed.
\end{examp}

\bigskip

Now, we are concerned with the sum rules for subdifferentials. For this, we
use the scalarization procedure described in Proposition \ref{pr_scal_subd}
and then we consider the real-valued functions associated to real-valued
set-valued maps, as observed in the opening of this section. Three preliminary
remarks, together with a lemma, are in order.

\begin{rmk}
Let $F_{1},F_{2}:X\rightrightarrows Y$ be set-valued maps. Clearly, for every
$\overline{x}\in X$,%
\[
\widehat{\partial}F_{1}\left(  \overline{x}\right)  +\widehat{\partial}%
F_{2}\left(  \overline{x}\right)  \subset\widehat{\partial}\left(  F_{1}%
+F_{2}\right)  \left(  \overline{x}\right)  .
\]

\end{rmk}

\begin{rmk}
If $A,B\subset\mathbb{R}$ are nonempty, closed and bounded from below, then
$A\subset B+[0,\infty)$ if and only if $\min A\geq\min B.$
\end{rmk}

\begin{rmk}
If $f:X\rightarrow\mathbb{R}$ is a function and $F:X\rightrightarrows
\mathbb{R}$ is a set-valued map with closed and bounded from below values,
then it is easy to see that the following equivalence holds: $f\left(
x\right)  =\min F\left(  x\right)  ,$ for every $x\in X$ if and only if
$F\left(  x\right)  +[0,\infty)=f\left(  x\right)  +[0,\infty),$ for every
$x\in X.$
\end{rmk}

\begin{lm}
\label{lm_iscusc}Let $f:X\rightarrow\mathbb{R}$ be a function and
$F:X\rightrightarrows\mathbb{R}$ be a set-valued map such that
\[
F\left(  x\right)  +[0,\infty)=f\left(  x\right)  +[0,\infty),\text{ }\forall
x\in X.
\]
Let $\overline{x}\in X.$

(i) If $F$ is u.c. at $\overline{x},$ then $f$ is lower semicontinuous (in
short, l.s.c.) at $\overline{x}$.

(ii) If $F$ is l.c. at $\overline{x},$ then $f$ is upper semicontinuous (in
short, u.s.c.) at $\overline{x}$.
\end{lm}

\noindent\textbf{Proof. }(i) Suppose that $F$ is u.c. at $\overline{x}$ and
take $\alpha\in\mathbb{R}$ with $\alpha<f\left(  \overline{x}\right)  ,$ i.e.,
$f\left(  \overline{x}\right)  \in\left(  \alpha,\infty\right)  .$ Then, since
$F\left(  \overline{x}\right)  \subset F\left(  \overline{x}\right)
+[0,\infty)=f\left(  \overline{x}\right)  +[0,\infty),$ one has that $F\left(
\overline{x}\right)  \subset\left(  \alpha,\infty\right)  ,$ whence, by using
the upper continuity of $F,$ one gets that there is a neighborhood $U$ of
$\overline{x}$ such that for all $u\in U$ one has $F\left(  u\right)
\subset\left(  \alpha,\infty\right)  .$ Now, since $f\left(  x\right)  \in
F\left(  x\right)  +[0,\infty),$ for every $x\in X,$ one obtains that for all
$u\in U,$ $f\left(  u\right)  \in\left(  \alpha,\infty\right)  ,$ i.e.,
$\alpha<f\left(  u\right)  ,$ so $f$ is l.s.c. at $\overline{x}.$

(ii) Suppose that $F$ is l.c. at $\overline{x}$ and take $\alpha\in\mathbb{R}$
with $f\left(  \overline{x}\right)  <\alpha.$ According to the hypothesis,
$f\left(  \overline{x}\right)  \in F\left(  \overline{x}\right)  +[0,\infty),$
so there exists $\overline{y}\in F\left(  \overline{x}\right)  $ such that
$f\left(  \overline{x}\right)  \geq\overline{y},$ therefore, $\overline{y}\in
F\left(  \overline{x}\right)  \cap\left(  -\infty,\alpha\right)  .$ Further,
by using the lower continuity of $F,$ one gets that there is a neighborhood
$U$ of $\overline{x}$ such that for all $u\in U$ one has $F\left(  u\right)
\cap\left(  -\infty,\alpha\right)  \neq\emptyset,$ so for all $u\in U,$ there
exists $v_{u}\in\mathbb{R}$ such that $v_{u}<\alpha$ and $v_{u}\in F\left(
u\right)  .$ Then, since $F\left(  x\right)  \subset f\left(  x\right)
+[0,\infty),$ for every $x\in X$ one has that $v_{u}\geq f\left(  u\right)  ,$
for all $u\in U,$ whence $f\left(  u\right)  <\alpha,$ for all $u\in U,$ and
the conclusion follows.\hfill$\square$

\bigskip

Notice that the reverse implications from the above lemma, in general, do not
hold as we illustrate in the below example.

\begin{examp}
Let $F_{1},F_{2}:\mathbb{R}\rightrightarrows\mathbb{R}$ be two set-valued maps
given by%
\[
F_{1}\left(  x\right)  =\left\{
\begin{array}
[c]{l}%
\left\{  -1,1\right\}  ,\text{ if }x\neq0,\\
\left\{  0\right\}  ,\text{ if }x=0,
\end{array}
\right.
\]
and%
\[
F_{2}\left(  x\right)  =\left\{
\begin{array}
[c]{l}%
\left\{  -1,-\frac{1}{2}\right\}  ,\text{ if }x=0,\\
\left\{  0\right\}  ,\text{ if }x\neq0,
\end{array}
\right.
\]
and $f_{1},f_{2}:\mathbb{R}\rightarrow\mathbb{R}$ two functions defined by%
\[
f_{1}\left(  x\right)  =\left\{
\begin{array}
[c]{l}%
-1,\text{ if }x\neq0,\\
0,\text{ if }x=0,
\end{array}
\right.
\]
and%
\[
f_{2}\left(  x\right)  =\left\{
\begin{array}
[c]{l}%
-1,\text{ if }x=0,\\
0,\text{ if }x\neq0.
\end{array}
\right.
\]
It is easy to see that for $i\in\left\{  1,2\right\}  ,F_{i}\left(  x\right)
+[0,\infty)=f_{i}\left(  x\right)  +[0,\infty),$ for every $x\in X,$ $f_{1}$
is u.s.c. at $0,$ but $F_{1}$ is not l.c. at $0$ and $f_{2}$ is l.s.c. at $0,$
but $F_{2}$ is not u.c. at $0.$
\end{examp}

Now we are able to present a sum rule for the convex case.

\begin{pr}
\label{pr_sum_conv}Let $F_{1},F_{2}:X\rightrightarrows Y$ be upper $K-$convex
set-valued maps with $K-$sequentially compact values. Suppose that $F_{1}$ (or
$F_{2}$) is continuous at some $x\in X$ or $X$ is finite dimensional. Then,
for every $\overline{x}\in X$ and $y^{\ast}\in K^{+}\setminus\left\{
0\right\}  ,$%
\[
y^{\ast}\left(  \widehat{\partial}\left(  F_{1}+F_{2}\right)  \left(
\overline{x}\right)  \right)  \subset\widehat{\partial}\left(  y^{\ast}\circ
F_{1}\right)  \left(  \overline{x}\right)  +\widehat{\partial}\left(  y^{\ast
}\circ F_{2}\right)  \left(  \overline{x}\right)  .
\]

\end{pr}

\noindent\textbf{Proof. }Take $T\in\widehat{\partial}\left(  F_{1}%
+F_{2}\right)  \left(  \overline{x}\right)  =\widehat{\partial}\left(
\operatorname*{Epi}F_{1}+\operatorname*{Epi}F_{2}\right)  \left(  \overline
{x}\right)  .$ According to Proposition \ref{pr_scal_subd}, for all $y^{\ast
}\in K^{+}\setminus\left\{  0\right\}  ,$%
\[
y^{\ast}\circ T\in\widehat{\partial}\left(  y^{\ast}\circ\left(
\operatorname*{Epi}F_{1}+\operatorname*{Epi}F_{2}\right)  \right)  \left(
\overline{x}\right)  =\widehat{\partial}\left(  y^{\ast}\circ
\operatorname*{Epi}F_{1}+y^{\ast}\circ\operatorname*{Epi}F_{2}\right)  \left(
\overline{x}\right)  .
\]
For $i\in\left\{  1,2\right\}  ,$ one has that $y^{\ast}\circ
\operatorname*{Epi}F_{i}$ has as values closed and bounded from below
intervals of $\mathbb{R}.$ Indeed, $\operatorname*{Epi}F_{i}$ has convex
values, whence $y^{\ast}\circ\operatorname*{Epi}F_{i}$ has as values intervals
in $\mathbb{R}$. Now, the $K-$sequentially compactness of the values of
$F_{i}$ implies their $K-$boundedness and for all $x\in X\ $and $i\in\left\{
1,2\right\}  $ we denote by $M_{x}^{i}$ the bounded set with the property
$F_{i}\left(  x\right)  \subset M_{x}^{i}+K.$ Then one has that for all
$z\in\left(  y^{\ast}\circ\operatorname*{Epi}F_{i}\right)  \left(  x\right)
,$ there is $y\in M_{x}^{i}$ and $k\in K$ such that $z=y^{\ast}\left(
y\right)  +y^{\ast}\left(  k\right)  ,$ whence%
\[
z\geq-\left\Vert y^{\ast}\right\Vert \left\Vert y\right\Vert ,
\]
and since $M_{x}^{i}$ is bounded, one deduces that $\left(  y^{\ast}%
\circ\operatorname*{Epi}F_{i}\right)  \left(  x\right)  $ are bounded from
below intervals. Finally, we show that these sets are also closed. Take again
$x\in X$ and a sequence $\left(  z_{n}\right)  \subset\left(  y^{\ast}%
\circ\operatorname*{Epi}F_{i}\right)  \left(  x\right)  $ such that
$z_{n}\rightarrow z.$ As above, for all $n$ there is $y_{n}\in F_{i}\left(
x\right)  $ and $k_{n}\in K$ such that $z_{n}=y^{\ast}\left(  y_{n}\right)
+y^{\ast}\left(  k_{n}\right)  .$ By the $K-$sequentially compactness of
$F_{i}\left(  x\right)  ,$ there is a sequence $\left(  c_{n}\right)  \subset
K$ such that $\left(  y_{n}-c_{n}\right)  $ converges towards some $u\in
F_{i}\left(  x\right)  .$ Therefore,%
\[
z_{n}=y^{\ast}\left(  y_{n}\right)  +y^{\ast}\left(  k_{n}\right)  =y^{\ast
}\left(  y_{n}-c_{n}\right)  +y^{\ast}\left(  c_{n}\right)  +y^{\ast}\left(
k_{n}\right)  ,
\]
so%
\[
z_{n}-y^{\ast}\left(  y_{n}-c_{n}\right)  \rightarrow z-y^{\ast}\left(
u\right)  \geq0.
\]
We get from here that $z\in y^{\ast}\left(  u\right)  +[0,\infty
)\subset\left(  y^{\ast}\circ F_{i}\right)  \left(  x\right)  +[0,\infty
)=\left(  y^{\ast}\circ\operatorname*{Epi}F_{i}\right)  \left(  x\right)  ,$
and the claim follows.

Observe now that the values of $y^{\ast}\circ\left(  \operatorname*{Epi}%
F_{1}+\operatorname*{Epi}F_{2}\right)  $ share the same properties and take
the minimal functions associated to these three set-valued maps, denoted as
$f_{F_{1}+F_{2}},$ $f_{F_{1}},$ and $f_{F_{2}},$ respectively. Observe as well
that these functions are convex, $f_{F_{1}+F_{2}}=f_{F_{1}}+f_{F_{2}}$ and at
any $x\in X,$ the subdifferential of $y^{\ast}\circ\left(  \operatorname*{Epi}%
F_{1}+\operatorname*{Epi}F_{2}\right)  $ coincides with the subdifferential of
$f_{F_{1}+F_{2}}$ and similarly for the other two functions. Then we have that%
\[
y^{\ast}\circ T\in\widehat{\partial}\left(  f_{F_{1}+F_{2}}\right)  \left(
\overline{x}\right)  =\widehat{\partial}\left(  f_{F_{1}}+f_{F_{2}}\right)
\left(  \overline{x}\right)  .
\]
If $X$ is finite dimensional, since $f_{_{F_{1}}}$ and $f_{F_{2}}$ are convex,
both functions $f_{_{F_{1}}}$ and $f_{F_{2}}$ are continuous. If $F_{1}$ is
continuous at $x\in X,$ then $y^{\ast}\circ F_{1}$ is continuous at $x$,
whence, by Lemma \ref{lm_iscusc}, $f_{_{F_{1}}}$ is continuous at $x$. This
means that the classical sum rule for the convex subdifferential of
real-valued convex functions can be applied, whence we get that for every
$\overline{x}\in X,$%
\[
y^{\ast}\circ T\in\widehat{\partial}f_{F_{1}}\left(  \overline{x}\right)
+\widehat{\partial}f_{F_{2}}\left(  \overline{x}\right)  =\widehat{\partial
}\left(  y^{\ast}\circ F_{1}\right)  \left(  \overline{x}\right)
+\widehat{\partial}\left(  y^{\ast}\circ F_{2}\right)  \left(  \overline
{x}\right)  .
\]
The conclusion follows.\hfill$\square$

\bigskip

As it is well known, the Fr\'{e}chet subdifferential does not enjoy exact sum
rules for general functions, so it is not expected to behave differently in
this more general setting. For this reason, we introduce now a limiting
subdifferential of a set-valued map in a similar manner as in the classical case.

Let $X$ be a Banach space. For $F:X\rightrightarrows Y$ with nonempty values
and $\overline{x}\in X$ define%
\[
\partial F\left(  \overline{x}\right)  =\left\{  T\in B\left(  X,Y\right)
\mid\exists\left(  x_{n}\right)  \overset{\operatorname*{Epi}F,e}{\rightarrow
}\overline{x},\text{ }\exists\left(  T_{n}\right)  \overset{SOT}{\rightarrow
}T,\text{ }T_{n}\in\widehat{\partial}F\left(  x_{n}\right)  ,\text{ }\forall
n\right\}  ,
\]
where $\left(  x_{n}\right)  \overset{\operatorname*{Epi}F,e}{\rightarrow
}\overline{x}$ means that $\left(  x_{n}\right)  \rightarrow\overline{x}$ and
$e\left(  \operatorname*{Epi}F\left(  x_{n}\right)  ,\operatorname*{Epi}%
F\left(  \overline{x}\right)  \right)  \rightarrow0,$ and $T_{n}\overset
{SOT}{\rightarrow}T$ means the convergence in the strong operator topology
(that is, $T_{n}\left(  x\right)  \rightarrow T\left(  x\right)  ,$ for all
$x\in X$).

This concept is a variation of the construction of the classical Mordukhovich
subdifferential in the context of Asplund spaces, but in a slightly weaker
form, since both the requirements $e\left(  \operatorname*{Epi}F\left(
x_{n}\right)  ,\operatorname*{Epi}F\left(  \overline{x}\right)  \right)
\rightarrow0$ and $e\left(  \operatorname*{Epi}F\left(  \overline{x}\right)
,\operatorname*{Epi}F\left(  x_{n}\right)  \right)  \rightarrow0$ would lead
in the case of real-valued functions to
\[
\partial f\left(  \overline{x}\right)  =\left\{  x^{\ast}\in X^{\ast}%
\mid\exists\left(  x_{n}\right)  \rightarrow\overline{x},\text{ }f\left(
x_{n}\right)  \rightarrow f\left(  \overline{x}\right)  ,\text{ }%
\exists\left(  x_{n}^{\ast}\right)  \overset{w^{\ast}}{\rightarrow}x^{\ast
},\text{ }x_{n}^{\ast}\in\widehat{\partial}f\left(  x_{n}\right)  ,\text{
}\forall n\right\}  ,
\]
where by $w^{\ast}$ we mean the weak star topology on $X^{\ast}.$ For more
details, the reader is referred to \cite[Theorem 2.34]{Morduk2006}.

\begin{rmk}
\label{rmk_lim}(i) Clearly, for all $x$ one has $\widehat{\partial}F\left(
x\right)  \subset\partial F\left(  x\right)  $ and $\partial F\left(
x\right)  =\partial\left(  \operatorname*{Epi}F\right)  \left(  x\right)  .$

(ii) If $F$ is from $X$ to $\mathbb{R}$ and has closed, bounded from below
values, then taking $f:X\rightarrow\mathbb{R}$ given by $f\left(  x\right)
=\min F\left(  x\right)  ,$ one has, for every $\overline{x}\in X,$ that
$t\in\partial F\left(  \overline{x}\right)  $ if and only if $t\in\partial
f\left(  \overline{x}\right)  .$
\end{rmk}

Clearly, convergence in the strong operator topology implies the convergence
in the weak operator topology (that is, $\left(  y^{\ast}\circ T_{n}\right)  $
converges pointwise to $y^{\ast}\circ T,$ for every $y^{\ast}\in Y^{\ast}$)
and, of course, it reduces to $w^{\ast}$ convergence when $Y=\mathbb{R}$.

\begin{pr}
\label{pr_scal_subd_l}Let $X$ be a Banach space, $F:X\rightrightarrows Y$ a
set-valued map, $\overline{x}\in X,$ and $T\in\partial F\left(  \overline
{x}\right)  $. Then for all $y^{\ast}\in K^{+}\setminus\left\{  0\right\}  ,$
$y^{\ast}\circ T\in\partial\left(  y^{\ast}\circ F\right)  \left(
\overline{x}\right)  .$
\end{pr}

\noindent\textbf{Proof. }Since $T\in\partial F\left(  \overline{x}\right)  ,$
there exist $\left(  x_{n}\right)  \overset{\operatorname*{Epi}F,e}%
{\rightarrow}\overline{x},$ $\left(  T_{n}\right)  \overset{SOT}{\rightarrow
}T$ such that for all $n,$ $T_{n}\in\widehat{\partial}F\left(  x_{n}\right)
.$ Following Proposition \ref{pr_scal_subd}, for all $n$ and all $y^{\ast}\in
K^{+}\setminus\left\{  0\right\}  ,$ $y^{\ast}\circ T_{n}\in\widehat{\partial
}\left(  y^{\ast}\circ F\right)  \left(  x_{n}\right)  .$ Fix $y^{\ast}\in
K^{+}\setminus\left\{  0\right\}  .$ Since $\left(  x_{n}\right)
\overset{\operatorname*{Epi}F,e}{\rightarrow}\overline{x},$ then $\left(
x_{n}\right)  \overset{\operatorname*{Epi}\left(  y^{\ast}\circ F\right)
,e}{\rightarrow}\overline{x}$ and since $\left(  T_{n}\right)  \overset
{SOT}{\rightarrow}T$, $\left(  y^{\ast}\circ T_{n}\right)  \overset
{SOT}{\rightarrow}y^{\ast}\circ T.$ This proves that $y^{\ast}\circ
T\in\partial\left(  y^{\ast}\circ F\right)  \left(  \overline{x}\right)
.$\hfill$\square$

\begin{pr}
Suppose that $X$ is a Banach space, $F$ is upper $K-$convex with
$K-$sequentially compact values around $\overline{x}\in X$. Assume that
$D_{Y}$ is $K-$closed. Then, one has that
\[
\partial F\left(  \overline{x}\right)  =\left\{  T\in B\left(  X,Y\right)
\mid F\left(  x\right)  \subset F\left(  \overline{x}\right)  +T\left(
x-\overline{x}\right)  +K,\text{ }\forall x\in X\right\}  .
\]

\end{pr}

\noindent\textbf{Proof. }The inclusion from right to left follows from Remark
\ref{rmk_lim} and Proposition \ref{pr_subd_F_convex}. Take now $T\in\partial
F\left(  \overline{x}\right)  .$ Then there exist $\left(  x_{n}\right)
\overset{\operatorname*{Epi}F,e}{\rightarrow}\overline{x},$ $\left(
T_{n}\right)  \overset{SOT}{\rightarrow}T$ such that for all $n,$ $T_{n}%
\in\widehat{\partial}F\left(  x_{n}\right)  .$ We can suppose, without loss of
generality, that $F\left(  x_{n}\right)  $ is $K-$sequentially compact for all
$n.$ According to Proposition \ref{pr_subd_F_convex},
\[
F\left(  x\right)  \subset F\left(  x_{n}\right)  +T_{n}\left(  x-x_{n}%
\right)  +K,\text{ }\forall x\in X,\text{ }\forall n.
\]
Observe that $T_{n}\left(  x-x_{n}\right)  \rightarrow T\left(  x-\overline
{x}\right)  .$ Indeed, by the Uniform Boundedness Principle, $\left(
\left\Vert T_{n}\right\Vert \right)  $ is a bounded sequence (notice that $X$
is a Banach space). Then the inequality
\begin{align*}
\left\Vert T_{n}\left(  x-x_{n}\right)  -T\left(  x-\overline{x}\right)
\right\Vert  &  \leq\left\Vert \left(  T_{n}-T\right)  \left(  x\right)
\right\Vert +\left\Vert T_{n}\left(  x_{n}\right)  -T\left(  \overline
{x}\right)  \right\Vert \\
&  \leq\left\Vert \left(  T_{n}-T\right)  \left(  x\right)  \right\Vert
+\left\Vert T_{n}\left(  x_{n}-\overline{x}\right)  \right\Vert +\left\Vert
T_{n}\left(  \overline{x}\right)  -T\left(  \overline{x}\right)  \right\Vert
\\
&  \leq\left\Vert \left(  T_{n}-T\right)  \left(  x\right)  \right\Vert
+\left\Vert T_{n}\right\Vert \left\Vert x_{n}-\overline{x}\right\Vert
+\left\Vert T_{n}\left(  \overline{x}\right)  -T\left(  \overline{x}\right)
\right\Vert ,\text{ }\forall n
\end{align*}
and the facts that $\left(  x_{n}\right)  \rightarrow\overline{x}$ and
$T_{n}\overset{SOT}{\rightarrow}T$ prove the thesis.

Take $\varepsilon>0.$ Then for all $x\in X$ and for all $n$ large enough,%
\[
F\left(  x\right)  \subset F\left(  x_{n}\right)  +T_{n}\left(  x-x_{n}%
\right)  +K\subset F\left(  \overline{x}\right)  +T\left(  x-\overline
{x}\right)  +\varepsilon D_{Y}+K.
\]
Since $\varepsilon$ is arbitrary and $F\left(  \overline{x}\right)  +T\left(
x-\overline{x}\right)  +K$ is closed, we get the conclusion.\hfill$\square$

\begin{pr}
\label{pr_sum_lim}Let $X$ be an Asplund space and $F_{1},F_{2}%
:X\rightrightarrows Y$ be set-valued maps with $K-$sequentially compact
values. Take $\overline{x}\in X$. Suppose that $F_{1}$ is u.c. around
$\overline{x},$ and $F_{2}$ is $K-$Lipschitz around $\overline{x}.$ Then for
every $y^{\ast}\in K^{+}\setminus\left\{  0\right\}  ,$%
\[
y^{\ast}\left(  \partial\left(  F_{1}+F_{2}\right)  \left(  \overline
{x}\right)  \right)  \subset\partial\left(  y^{\ast}\circ F_{1}\right)
\left(  \overline{x}\right)  +\partial\left(  y^{\ast}\circ F_{2}\right)
\left(  \overline{x}\right)  .
\]

\end{pr}

\noindent\textbf{Proof. }Take $T\in\partial\left(  F_{1}+F_{2}\right)  \left(
\overline{x}\right)  =\partial\left(  \operatorname*{Epi}F_{1}%
+\operatorname*{Epi}F_{2}\right)  \left(  \overline{x}\right)  .$ By
Proposition \ref{pr_scal_subd_l}, for all $y^{\ast}\in K^{+}\setminus\left\{
0\right\}  ,$%
\[
y^{\ast}\circ T\in\partial\left(  y^{\ast}\circ\left(  \operatorname*{Epi}%
F_{1}+\operatorname*{Epi}F_{2}\right)  \right)  \left(  \overline{x}\right)
.
\]
As in the proof of Proposition \ref{pr_sum_conv}, $y^{\ast}\circ
\operatorname*{Epi}F_{i}$ ($i\in\left\{  1,2\right\}  $) has closed and
bounded from below values (subsets of $\mathbb{R}$) and the same can be said
for the set-valued map $y^{\ast}\circ\left(  \operatorname*{Epi}%
F_{1}+\operatorname*{Epi}F_{2}\right)  .$ Then we can as well consider the
minimal functions associated to these three set-valued maps, denoted as in the
mentioned proof. According to Remark \ref{rmk_lim},
\[
y^{\ast}\circ T\in\partial\left(  f_{F_{1}+F_{2}}\right)  \left(  \overline
{x}\right)  =\partial\left(  f_{F_{1}}+f_{F_{2}}\right)  \left(  \overline
{x}\right)  .
\]
Again by Lemma \ref{lm_iscusc}, since $F_{1}$ is u.c. around $\overline{x},$
we get that $f_{F_{1}}$ is l.s.c. around this point. Similarly, since $F_{2}$
is $K-$Lipschitz around $\overline{x},$ $f_{F_{2}}$ is Lipschitz around this
point. Then one can apply the sum rule for the limiting subdifferential (see
\cite[Theorem 2.33]{Morduk2006}) to get%
\[
y^{\ast}\circ T\in\partial f_{F_{1}}\left(  \overline{x}\right)  +\partial
f_{F_{2}}\left(  \overline{x}\right)  ,
\]
and the conclusion ensues.\hfill$\square$

\bigskip

Finally, we are interested in some statements concerning the subdifferential
of some special (multi)functions.

For $e\in K\setminus\left\{  0\right\}  $ we use the notation%
\[
K_{e}^{+}=\left\{  y^{\ast}\in K^{+}\mid y^{\ast}\left(  e\right)
\neq0\right\}  .
\]
Notice that $K_{e}^{+}\neq\emptyset$ (by a standard separation theorem applied
to $\left\{  e\right\}  $ and $-K$) and, obviously, if $\operatorname*{int}%
K\neq\emptyset$ and $e\in\operatorname*{int}K,$ then $K_{e}^{+}=K^{+}%
\setminus\left\{  0\right\}  .$

\begin{rmk}
\label{rmk_f_e}Take in Proposition \ref{pr_scal_subd} the mapping $F:=f$ as a
function of the form $\varphi\left(  \cdot\right)  e,$ where $\varphi
:X\rightarrow\mathbb{R}$ and $e\in K\setminus\left\{  0\right\}  .$ Therefore,
if $t\in\widehat{\partial}f\left(  \overline{x}\right)  $, then for all
$y^{\ast}\in K_{e}^{+}$ one has that $y^{\ast}\left(  e\right)  ^{-1}\left(
y^{\ast}\circ t\right)  \in\widehat{\partial}\varphi\left(  \overline
{x}\right)  $.
\end{rmk}

\begin{pr}
\label{pr_nevid}Let $\varphi:X\rightarrow\mathbb{R}$ be a function,
$\overline{x}\in X,$ and $e\in K\setminus\left\{  0\right\}  .$ Denote
$f=\varphi\left(  \cdot\right)  e.$ Then,
\[
\left\{  x^{\ast}\left(  \cdot\right)  e\mid x^{\ast}\in\widehat{\partial
}\varphi\left(  \overline{x}\right)  \right\}  \subset\widehat{\partial
}f\left(  \overline{x}\right)  .
\]

\end{pr}

\noindent\textbf{Proof. }Let $x^{\ast}\in\widehat{\partial}\varphi\left(
\overline{x}\right)  .$ Then for all $\varepsilon>0,$ there is $\delta>0$ such
that for all $x\in B\left(  \overline{x},\delta\right)  ,$
\[
\varphi\left(  x\right)  \in\varphi\left(  \overline{x}\right)  +x^{\ast
}\left(  x-\overline{x}\right)  +\varepsilon\left\Vert x-\overline
{x}\right\Vert \left[  -1,1\right]  +[0,\infty).
\]
Then, for such $x,$
\begin{align*}
\varphi\left(  x\right)  e  &  \in\varphi\left(  \overline{x}\right)
e+x^{\ast}\left(  x-\overline{x}\right)  e+\varepsilon\left\Vert
x-\overline{x}\right\Vert \left[  -1,1\right]  e+[0,\infty)e\\
&  \subset\varphi\left(  \overline{x}\right)  e+x^{\ast}\left(  x-\overline
{x}\right)  e+\varepsilon\left\Vert x-\overline{x}\right\Vert \left\Vert
e\right\Vert D_{Y}+K,
\end{align*}
and this is enough to prove that $x^{\ast}\left(  \cdot\right)  e\in
\widehat{\partial}f\left(  \overline{x}\right)  .$\hfill$\square$

\begin{pr}
\label{pr_calc_ne}Let $e\in K\setminus\left\{  0\right\}  \ $and denote
$f=\left\Vert \cdot\right\Vert e$ (that is, $\varphi=\left\Vert \cdot
\right\Vert $ with the above notation). Then, $D_{X^{\ast}}e\subset
\widehat{\partial}f\left(  0\right)  $ and for all $y^{\ast}\in K_{e}^{+},$
$y^{\ast}\left(  \widehat{\partial}f\left(  0\right)  \right)  \subset
y^{\ast}\left(  e\right)  D_{X^{\ast}}.$
\end{pr}

\noindent\textbf{Proof. }By Proposition \ref{pr_nevid}, $D_{X^{\ast}}%
e\subset\widehat{\partial}f\left(  0\right)  $, whence the first inclusion
holds. Take $t\in\widehat{\partial}f\left(  0\right)  .$ By Remark
\ref{rmk_f_e}, for all $y^{\ast}\in K_{e}^{+}$ one has $y^{\ast}\left(
e\right)  ^{-1}\left(  y^{\ast}\circ t\right)  \in\widehat{\partial}%
\varphi\left(  0\right)  ,$ so there is $x_{y^{\ast}}^{\ast}\in D_{X^{\ast}}$
such that $y^{\ast}\left(  e\right)  ^{-1}\left(  y^{\ast}\circ t\right)
=x_{y^{\ast}}^{\ast},$ which implies the conclusion.\hfill$\square$

\bigskip

Finally, inspired by \cite{MNY}, we prove a difference rule.

\begin{pr}
\label{pr_diff}Let $F:X\rightrightarrows Y$ be a set-valued map,
$\varphi:X\rightarrow\mathbb{R}$ a function, $\overline{x}\in X,$ and $e\in
K\setminus\left\{  0\right\}  .$ Denote $f=\varphi\left(  \cdot\right)  e.$ If
$\widehat{\partial}\varphi\left(  \overline{x}\right)  \neq\emptyset,$ then
\[
\widehat{\partial}\left(  F-f\right)  \left(  \overline{x}\right)
\subset\bigcap_{\substack{t\in\widehat{\partial}f\left(  \overline{x}\right)
\\y^{\ast}\in K_{e}^{+}}}\left(  \widehat{\partial}F\left(  \overline
{x}\right)  -y^{\ast}\left(  e\right)  ^{-1}\left(  y^{\ast}\circ t\right)
e\right)  .
\]

\end{pr}

\noindent\textbf{Proof. }Take $T\in\widehat{\partial}\left(  F-f\right)
\left(  \overline{x}\right)  .$ Then, for all $\varepsilon>0$ there is
$\delta>0$ such that for all$\ x\in B\left(  \overline{x},\delta\right)  $ one
has%
\[
F\left(  x\right)  -\varphi\left(  x\right)  e\subset F\left(  \overline
{x}\right)  -\varphi\left(  \overline{x}\right)  e+T\left(  x-\overline
{x}\right)  +\varepsilon\left\Vert x-\overline{x}\right\Vert D_{Y}+K.
\]
Since $\widehat{\partial}\varphi\left(  \overline{x}\right)  \neq\emptyset,$
by Proposition \ref{pr_nevid}, $\widehat{\partial}f\left(  \overline
{x}\right)  \neq\emptyset.$ Take $t\in\widehat{\partial}f\left(  \overline
{x}\right)  $ and $y^{\ast}\in K_{e}^{+}.$ According to Remark \ref{rmk_f_e}%
$,$%
\[
\left(  y^{\ast}\left(  e\right)  \right)  ^{-1}\left(  y^{\ast}\circ
t\right)  \in\widehat{\partial}\varphi\left(  \overline{x}\right)  .
\]
Using the variational description of the Fr\'{e}chet subgradients, and making
$\delta$ smaller if necessary, there is a function $s_{y^{\ast}},$
differentiable at $\overline{x},$ such that $s_{y^{\ast}}\left(  \overline
{x}\right)  =\varphi\left(  \overline{x}\right)  ,$ $\nabla s_{y^{\ast}%
}\left(  \overline{x}\right)  =y^{\ast}\left(  e\right)  ^{-1}\left(  y^{\ast
}\circ t\right)  $ and $s_{y^{\ast}}\left(  x\right)  \leq\varphi\left(
x\right)  $ for all $x\in B\left(  \overline{x},\delta\right)  .$ By the
differentiability of $s_{y^{\ast}}$ at $\overline{x},$ there is a function
$\alpha:X\rightarrow\mathbb{R}$ continuous at $\overline{x}$ and with
$\lim_{x\rightarrow\overline{x}}\alpha\left(  x\right)  =0$ such that for all
$x,$%
\[
s_{y^{\ast}}\left(  x\right)  -s_{y^{\ast}}\left(  \overline{x}\right)
=\nabla s_{y^{\ast}}\left(  \overline{x}\right)  \left(  x-\overline
{x}\right)  +\alpha\left(  x\right)  \left\Vert x-\overline{x}\right\Vert .
\]
We can suppose, again making $\delta$ smaller if necessary, that
$\alpha\left(  x\right)  e\in\varepsilon D_{Y}$ for all $x\in B\left(
\overline{x},\delta\right)  .$ We successively get, for all $x\in B\left(
\overline{x},\delta\right)  ,$%
\begin{align*}
F\left(  x\right)   &  \subset F\left(  \overline{x}\right)  +\left(
\varphi\left(  x\right)  -\varphi\left(  \overline{x}\right)  \right)
e+T\left(  x-\overline{x}\right)  +\varepsilon\left\Vert x-\overline
{x}\right\Vert D_{Y}+K\\
&  \subset F\left(  \overline{x}\right)  +\left(  s_{y^{\ast}}\left(
x\right)  -s_{y^{\ast}}\left(  \overline{x}\right)  \right)  e+T\left(
x-\overline{x}\right)  +\varepsilon\left\Vert x-\overline{x}\right\Vert
D_{Y}+K\\
&  \subset F\left(  \overline{x}\right)  +\left(  \nabla s_{y^{\ast}}\left(
\overline{x}\right)  \left(  x-\overline{x}\right)  +\alpha\left(  x\right)
\left\Vert x-\overline{x}\right\Vert \right)  e+T\left(  x-\overline
{x}\right)  +\varepsilon\left\Vert x-\overline{x}\right\Vert D_{Y}+K\\
&  =F\left(  \overline{x}\right)  +\left(  T\left(  \cdot\right)  +\nabla
s_{y^{\ast}}\left(  \overline{x}\right)  \left(  \cdot\right)  e\right)
\left(  x-\overline{x}\right)  +\varepsilon\left\Vert x-\overline
{x}\right\Vert D_{Y}+\alpha\left(  x\right)  \left\Vert x-\overline
{x}\right\Vert e+K\\
&  \subset F\left(  \overline{x}\right)  +\left(  T\left(  \cdot\right)
+\nabla s_{y^{\ast}}\left(  \overline{x}\right)  \left(  \cdot\right)
e\right)  \left(  x-\overline{x}\right)  +2\varepsilon\left\Vert
x-\overline{x}\right\Vert D_{Y}+K.
\end{align*}
From this, we get that
\[
T+\nabla s_{y^{\ast}}\left(  \overline{x}\right)  e\in\widehat{\partial
}F\left(  \overline{x}\right)  ,
\]
whence
\[
T\in\widehat{\partial}F\left(  \overline{x}\right)  -\nabla s_{y^{\ast}%
}\left(  \overline{x}\right)  e=\widehat{\partial}F\left(  \overline
{x}\right)  -y^{\ast}\left(  e\right)  ^{-1}\left(  y^{\ast}\circ t\right)
e,
\]
and the conclusion follows.\hfill$\square$

\section{Ideal solutions}

Many authors have discussed optimization problems that are governed by
multifunctions, based on well-established relations on sets defined by Kuroiwa
(see \cite{KN}). We recall one such relation here and a concept of minimality
associated with it (see \cite{BDS}). For two nonempty sets
\[
A\preceq_{K}^{l}B\iff B\subset A+K.
\]

\begin{df}
\label{def_min}Let $F:X\rightrightarrows Y$ be a set valued-map with nonempty
values and $M\subset X$ be a nonempty set. One says that $\overline{x}\in M$
is a local $l-$minimum for $F$ on $M$ if there exists $\varepsilon>0$ such
that%
\[
x\in M\cap B\left(  \overline{x},\varepsilon\right)  ,\text{ }F\left(
x\right)  \preceq_{K}^{l}F\left(  \overline{x}\right)  \Longrightarrow
F\left(  \overline{x}\right)  \preceq_{K}^{l}F\left(  x\right)  .
\]

\end{df}

\begin{examp}
Notice that, in general, the local $\ell-$minimality of $\overline{x}$ for $F$
does not imply $0\in\widehat{\partial}F\left(  \overline{x}\right)  .$ For
example, $\overline{x}=0$ is a $l-$minimum for $F$ on $M=\mathbb{R},$ but
$0\notin\widehat{\partial}F\left(  \overline{x}\right)  ,$ where
$K:=\mathbb{R}_{+}^{2},$ $F:\mathbb{R}\rightrightarrows\mathbb{R}^{2},$%
\[
F\left(  x\right)  =\left\{
\begin{array}
[c]{l}%
\left\{  x\right\}  \times\lbrack-2,\infty),\text{ if }x>0,\\
\left\{  0\right\}  \times\lbrack0,\infty),\text{ if }x\leq0.
\end{array}
\right.
\]
However, if $F$ is from $X$ to $\mathbb{R}$ and has closed, bounded from below
values, then taking $f:X\rightarrow\mathbb{R}$ given by $f\left(  x\right)
=\min F\left(  x\right)  ,$ one has that $\ell-$minimality of $\overline{x}$
for $F$ is equivalent to classical minimality of $\overline{x}$ for $f,$ so in
this case $0\in\widehat{\partial}F\left(  \overline{x}\right)  .$
\end{examp}

We consider now a concept for which the basic Fermat rule takes place.

\begin{df}
Let $M\subset X$ be a nonempty set. One says that $\overline{x}\in M$ is a
local ideal minimum for $F$ on $M$ if there is $\varepsilon>0$ such that%
\begin{equation}
F\left(  \overline{x}\right)  \not \subset F\left(  x\right)  +Y\setminus
-K,\text{ }\forall x\in M\cap B\left(  \overline{x},\varepsilon\right)
\setminus\left\{  \overline{x}\right\}  . \label{def_id}%
\end{equation}
The global version of the above concept is obtained by taking the open ball
$B\left(  \overline{x},\varepsilon\right)  $ as the whole space $X,$ in which
case we omit to write "local". Also, we omit to write "on $M$" when we study
the case without restrictions.
\end{df}

\begin{rmk}
\label{0apartine_dif}Notice that relation (\ref{def_id}) implies%
\[
F\left(  x\right)  \subset F\left(  \overline{x}\right)  +K,\text{ }\forall
x\in M\cap B\left(  \overline{x},\varepsilon\right)  .
\]
Consequently, if $\overline{x}\in\operatorname*{int}M$ is a local ideal
minimum for $F,$ then $0\in\widehat{\partial}F\left(  \overline{x}\right)  .$
If $F:=f$ is a function, then the above minimality becomes%
\[
f\left(  x\right)  -f\left(  \overline{x}\right)  \in K,\text{ }\forall x\in
M\cap B\left(  \overline{x},\varepsilon\right)  ,
\]
and this is known in vector optimization under the name of ideal minimality.
\end{rmk}

\begin{lm}
\label{lm_comp}Let $A\subset Y\ $be a $K-$sequentially compact set and
$D\subset Y$ be a set such that $A\subset D+Y\setminus-K.$ Then there exists
$\varepsilon>0$ such that $A+B\left(  0,\varepsilon\right)  \subset
D+Y\setminus-K.$
\end{lm}

\noindent\textbf{Proof. }Suppose, by way of contradiction, that the conclusion
is not true. Then one can find a sequence $\left(  \rho_{n}\right)
\rightarrow0$ and $\left(  a_{n}\right)  \subset A$ such that for all
$n\geq1,$%
\[
a_{n}+\rho_{n}\notin D+Y\setminus-K.
\]
Since $A$ is $K-$sequentially compact, there is a sequence $\left(
c_{n}\right)  \subset K$ such that a subsequence $\left(  a_{n_{k}}-c_{n_{k}%
}\right)  $ converges towards some $a\in A.$ Then $a$ belongs to the open set
$D+Y\setminus-K,$ so for $k$ large enough,%
\[
a_{n_{k}}-c_{n_{k}}+\rho_{n_{k}}\in D+Y\setminus-K,
\]
which, having in view that $Y\setminus-K+K=Y\setminus-K,$ implies that%
\[
a_{n_{k}}+\rho_{n_{k}}\in D+Y\setminus-K,
\]
for large $k.$ This is a contradiction, whence the conclusion holds.\hfill
$\square$

\bigskip

In the following, we present a result which follow a well-established approach
known as penalization (as described by Clarke in \cite{clarke} and extended to
vector optimization in \cite{Ye}, and also to set optimization, e.g., in
\cite{BDS}). This approach involves incorporating the distance between the
constraint set and the objective map to convert geometric constraints into an
unconstrained problem. The Lipschitz assumption used in the below penalization
result is consistent with the one employed in previous studies of penalization
(see, e.g., \cite{Ye}, \cite{BDS}) and also it is in the line of the Lipschitz
notion given in Definition \ref{def_K_lip}.

\begin{pr}
\label{pr_pen_local}Let $M\subset X$ be a nonempty set, assume that
$\overline{x}\in M\cap M^{\prime},$ $F\left(  \overline{x}\right)  $ is
$K-$sequentially compact and $\overline{x}$ is a local ideal minimum for $F$
on $M.$ Suppose that the following generalized Lipschitz condition holds:
there are $e\in K\mathbb{\setminus}\left\{  0\right\}  ,$ $r>0,$ $\ell>0$ such
that for all $u\in\left(  X\setminus M\right)  \cap B\left(  \overline
{x},r\right)  $ and $v\in M\cap B\left(  \overline{x},r\right)  ,$%
\[
F\left(  u\right)  +\ell\left\Vert u-v\right\Vert e\subset F\left(  v\right)
+K.
\]
Then $\overline{x}$ is a local ideal minimum on $X$ (that is, without
constraints) for the set-valued map $G:X\rightrightarrows Y,$%
\[
G\left(  x\right)  =F\left(  x\right)  +\ell d_{M}\left(  x\right)  e.
\]

\end{pr}

\noindent\textbf{Proof.} Without loss of generality, one can take $r$ as the
radius of the ball around $\overline{x}$ where the minimality condition holds.
Suppose, by way of contradiction, that there is $u\in B\left(  \overline
{x},3^{-1}r\right)  $ such that%
\[
F\left(  \overline{x}\right)  \subset F\left(  u\right)  +\ell d_{M}\left(
u\right)  e+Y\setminus-K.
\]
Clearly, $u\notin M.$ By virtue of Lemma \ref{lm_comp}, there is
$\varepsilon>0$ such that%
\[
F\left(  \overline{x}\right)  -\ell\varepsilon e\subset F\left(  u\right)
+\ell d_{M}\left(  u\right)  e+Y\setminus-K,
\]
that is%
\[
F\left(  \overline{x}\right)  \subset F\left(  u\right)  +\ell\left(
d_{M}\left(  u\right)  +\varepsilon\right)  e+Y\setminus-K.
\]
On the other hand, by tacking $\delta=\min\left\{  \varepsilon,3^{-1}%
r\right\}  ,$ there exists $a\in M$ such that
\[
\left\Vert u-a\right\Vert <d_{M}\left(  u\right)  +\delta.
\]
Since $\overline{x}\in M^{\prime}$, we can ensure that exists $a\in M$ such
that $a\neq\overline{x}$ and the above inequality holds$,$ so $a\in
M\setminus\left\{  \overline{x}\right\}  $ and
\begin{align*}
\left\Vert a-\overline{x}\right\Vert  &  \leq\left\Vert a-u\right\Vert
+\left\Vert u-\overline{x}\right\Vert <d_{M}\left(  u\right)  +\delta
+\left\Vert u-\overline{x}\right\Vert \\
&  \leq2\left\Vert u-\overline{x}\right\Vert +3^{-1}r<r,
\end{align*}
whence $a\in\left(  M\setminus\left\{  \overline{x}\right\}  \right)  \cap
B\left(  \overline{x},r\right)  $. We have%
\[
\left(  d_{M}\left(  u\right)  +\varepsilon-\left\Vert u-a\right\Vert \right)
e=\left(  d_{M}\left(  u\right)  +\delta-\left\Vert u-a\right\Vert \right)
e+\left(  \varepsilon-\delta\right)  e\in K\mathbb{\setminus}\left\{
0\right\}  ,
\]
and we get%
\[
\left(  d_{M}\left(  u\right)  +\varepsilon\right)  e\in\left\Vert
u-a\right\Vert e+K\mathbb{\setminus}\left\{  0\right\}  .
\]
Consequently,%
\begin{align*}
F\left(  \overline{x}\right)   &  \subset F\left(  u\right)  +\ell\left\Vert
u-a\right\Vert e+K\mathbb{\setminus}\left\{  0\right\}  +Y\setminus-K\\
&  \subset F\left(  u\right)  +\ell\left\Vert u-a\right\Vert e+Y\setminus-K.
\end{align*}
Using the generalized Lipschitz property we consider as hypothesis, we get%
\[
F\left(  \overline{x}\right)  \subset F\left(  a\right)  +Y\setminus-K,
\]
which is a contradiction. So, we can deduce that the conclusion is
true.\hfill$\square$

\bigskip

The global version reads as follows.

\begin{pr}
\label{pr_pen}Let $M\subset X$ be a nonempty set, assume that $\overline{x}\in
M\cap M^{\prime},$ $F\left(  \overline{x}\right)  $ is $K-$sequentially
compact and $\overline{x}$ is an ideal minimum for $F$ on $M.$ Suppose that
the following generalized Lipschitz condition holds: there are $e\in
K\mathbb{\setminus}\left\{  0\right\}  $ and $\ell>0$ such that for all $u\in
X\setminus M$ and $v\in M,$%
\[
F\left(  u\right)  +\ell\left\Vert u-v\right\Vert e\subset F\left(  v\right)
+K.
\]
Then $\overline{x}$ is an ideal minimum on $X$ (that is, without constraints)
for the set-valued map $G:X\rightrightarrows Y,$%
\[
G\left(  x\right)  =F\left(  x\right)  +\ell d_{M}\left(  x\right)  e.
\]

\end{pr}

Finally, putting together the local penalization result with some of the
subdifferential calculus rules developed before, we write necessary optimality
conditions for ideal sharp minimality (for more details, see the comment after
the proof).

\begin{pr}
Let $X$ be an Asplund space, $M\subset X$ be a nonempty set and $\mu>0.$
Suppose that the following generalized Lipschitz condition holds: there are
$e\in K\mathbb{\setminus}\left\{  0\right\}  ,$ $r>0,$ $\ell>0$ such that for
all $u\in\left(  X\setminus M\right)  \cap B\left(  \overline{x},r\right)  $
and $v\in M\cap B\left(  \overline{x},r\right)  ,$%
\begin{equation}
F\left(  u\right)  +\ell\left\Vert u-v\right\Vert e\subset F\left(  v\right)
+K. \label{K-lip}%
\end{equation}
Take $\overline{x}\in M\cap M^{\prime}$ as a local ideal minimum for $F\left(
\cdot\right)  -\mu\left\Vert \cdot-\overline{x}\right\Vert e$ on $M$. Suppose
that $F$ has $K-$sequentially compact values and it is u.c. around
$\overline{x}.$ Then for all $y^{\ast}\in K_{e}^{+},$%
\[
\mu y^{\ast}\left(  e\right)  D_{X^{\ast}}\subset\partial\left(  y^{\ast}\circ
F\right)  \left(  \overline{x}\right)  +\left(  \ell+\mu\right)  y^{\ast
}\left(  e\right)  \partial d_{M}\left(  \overline{x}\right)  .
\]

\end{pr}

\noindent\textbf{Proof. }From relation (\ref{K-lip}) one gets that for all
$u\in\left(  X\setminus M\right)  \cap B\left(  \overline{x},r\right)  $ and
$v\in M\cap B\left(  \overline{x},r\right)  ,$%
\[
F\left(  u\right)  -\mu\left\Vert u-\overline{x}\right\Vert e+\left(  \ell
+\mu\right)  \left\Vert u-v\right\Vert e\subset F\left(  v\right)
-\mu\left\Vert v-\overline{x}\right\Vert e+K.
\]
i.e., the Lipschitz assumption used in the local penalization result from
above holds for the set-valued map $x\rightrightarrows F\left(  x\right)
-\mu\left\Vert x-\overline{x}\right\Vert e.$ Therefore, according to
Proposition \ref{pr_pen_local}, $\overline{x}\in M$ is a local ideal minimum
for $F\left(  \cdot\right)  -\mu\left\Vert \cdot-\overline{x}\right\Vert
e+\left(  \ell+\mu\right)  d_{M}\left(  \cdot\right)  e$ on $X.$ By
Proposition \ref{pr_diff} and Remark \ref{0apartine_dif}, one has%
\begin{align*}
0  &  \in\widehat{\partial}\left(  F\left(  \cdot\right)  +\left(  \ell
+\mu\right)  d_{M}\left(  \cdot\right)  e-\mu\left\Vert \cdot-\overline
{x}\right\Vert e\right)  \left(  \overline{x}\right) \\
&  \subset\bigcap_{\substack{t\in\widehat{\partial}\left(  \mu\left\Vert
\cdot-\overline{x}\right\Vert e\right)  \left(  \overline{x}\right)
\\y^{\ast}\in K_{e}^{+}}}\left(  \widehat{\partial}\left(  F+\left(  \ell
+\mu\right)  d_{M}e\right)  \left(  \overline{x}\right)  -y^{\ast}\left(
e\right)  ^{-1}\left(  y^{\ast}\circ t\right)  e\right)  .
\end{align*}
Using now Proposition \ref{pr_sum_lim}, for all $t\in\widehat{\partial}\left(
\mu\left\Vert \cdot-\overline{x}\right\Vert e\right)  \left(  \overline
{x}\right)  $, $y^{\ast}\in K_{e}^{+},$ and $z^{\ast}\in K^{+}\setminus
\left\{  0\right\}  ,$
\begin{align*}
z^{\ast}\left(  y^{\ast}\left(  e\right)  ^{-1}\left(  y^{\ast}\circ t\right)
e\right)   &  \in z^{\ast}\left(  \widehat{\partial}\left(  F+\left(  \ell
+\mu\right)  d_{M}e\right)  \left(  \overline{x}\right)  \right)  \subset
z^{\ast}\left(  \partial\left(  F+\left(  \ell+\mu\right)  d_{M}e\right)
\left(  \overline{x}\right)  \right) \\
&  \subset\partial\left(  z^{\ast}\circ F\right)  \left(  \overline{x}\right)
+\left(  \ell+\mu\right)  z^{\ast}\left(  e\right)  \partial d_{M}\left(
\overline{x}\right)  .
\end{align*}
Taking $z^{\ast}=y^{\ast}$, and using the first inclusion given in Proposition
\ref{pr_calc_ne}, the conclusion follows.\hfill$\square$

\begin{rmk}
Observe that the minimality for $F\left(  \cdot\right)  -\mu\left\Vert
\cdot-\overline{x}\right\Vert e$ corresponds to the extension in the current
setting of the concept of sharp minimum (see \cite{W} for the scalar case and
\cite{FJ} for the vectorial one).
\end{rmk}

\bigskip

\noindent\textbf{Funding.} This work was supported by a grant of the Ministry
of Research, Innovation and Digitization, CNCS - UEFISCDI, project number
PN-III-P4-PCE-2021-0690, within PNCDI III.

\noindent\textbf{Data availability.} This manuscript has no associated data.

\noindent\textbf{Disclosure statement}. No potential conflict of interest was
reported by the authors.

\end{document}